\documentclass{article}
\usepackage{amsmath,amssymb,amsthm}
\newtheorem*{thm}{Theorem}

\newtheorem{prop}{Proposition}
\newtheorem{defin}{Definition}

\newtheorem*{rk}{Remarks}

\begin{document}
\author{Ma{\l}gorzata Stawiska\\
Department of Mathematics, Statistics \& Computer Science\\
University of Illinois at Chicago\\
851 S. Morgan St.\\
Chicago, IL 60607\\
\texttt{stawiska@math.uic.edu}}
\title{Attracting divisors on projective algebraic varieties}
\maketitle
\begin{abstract}
We obtain sufficient and necessary conditions (in terms of
positive singular metrics on an associated line bundle) for a
positive divisor $D$ on a projective algebraic variety $X$ to be
attracting for a holomorphic map $f:X \mapsto X$.
\end{abstract}

\section{Preliminaries}\label{s:prelim}

Attracting sets play an important part in the study of dynamical
systems. In recent years much attention has been devoted to
attracting periodic points for maps in several complex
variables,(\cite{FS}, \cite{Ga}, \cite{Ue} and other) and, to
somewhat lesser extent, to attracting hypersurfaces (\cite{BD},
\cite{BDM}, \cite{St1}, \cite{St2}). It would be interesting to
find unifying framework for dealing with attracting periodic
points on Riemann surfaces and hypersurfaces in higher dimensional
complex manifolds as well as to undertake study of algebraic
attracting sets of codimension greater than 1. We hope that
initial steps towards both goals can be made by considering
attracting divisors for holomorphic maps on projective algebraic
varieties, which are the subject of the present paper. We
characterize an attracting divisor $D$ for a holomorphic map
$f:X\mapsto X$ by behavior of a suitable positive singular
hermitian metric on the line bundle $[D]$ and of an associated
class of $\omega$-plurisubharmonic functions. In sections
\ref{s:prelim} and \ref{s:metric} we gather theoretical tools and
in Section \ref{s:attract} we prove our main result:\\
\begin{thm} The following are equivalent:\\
(i) $D$ is attracting for $f$;\\
(ii) there is a positive singular metric with weight $\psi$ on the
line bundle $[D]$ and a neighborhood $V$ of $A$ such that $\exists
0 < \beta < 1: \quad
f_*\psi \geq \beta \psi$ in $V$;\\
(iii)there is a positive singular metric with weight $\phi$ on the
line bundle $[D]$ such that $\exists 0 < \beta < 1: \quad
f_*(PSH(X,\omega)) \subset \beta PSH(X,\omega)$, where $\omega =
dd^c\phi$.
\end{thm}
Let $X$ be an $n$-dimensional complex algebraic variety (we will
make more assumptions about $X$ later). $\mathcal{O}_X$ will
denote the sheaf of germs of holomorphic functions on $X$ and
$\mathcal{M}$ the sheaf of germs of meromorphic functions on $X$.
For detailed definitions and more information we refer the reader
to \cite{GR},
\cite{GH}, \cite{Fu} or \cite{Ha}. \\
\begin{defin}(\cite{Fu}, B.4.1, B.4.2, B.4.5) A (Cartier) divisor $D$ on $X$
is defined by data $(U_\iota,f_\iota)$, where the $U_\iota$ form
an open covering of $X$ and the $f_\iota$ are non-zero meromorphic
functions, such that $f_\iota/f_\kappa$ is a unit (i.e., a
holomorphic, nowhere vanishing function) on
$U_{\iota\kappa}=U_\iota \cap U_\kappa$. The rational functions
$f_\iota$ are called local equations for $D$ (they are determined
up to multiplication by units on $U_\iota$). The support of $D$ is
the set of all points $x \in X$ such that a local equation for $D$
is not in $\mathcal{O}^*_{x,X}$ (i.e., it is not a unit). A
divisor $D$ is effective if local equations $f_\iota$ are sections
of $\mathcal{O}$ on $U_\iota$.
\end{defin}

Writing the local equation of $D$ on $U_\iota$ as $f_\iota =
a_\iota/b_\iota$, with $a_\iota,b_\iota$ holomorphic on $U_\iota$,
we see that the support of $D$ in $U_\iota$ consists of components
of the set $a_\iota b_\iota = 0$. Hence (see \cite{L}, II.5.3),
the support of $D \neq 0$ in $X$ is an analytic
subset of $X$ of pure codimension one.\\

We will use interchangeably the languages of divisors, line
bundles and invertible sheaves, because of the following one-to-one correspondences
(cf. \cite{Fu}, Appendix B.4.4 and \cite{Ha}, Ch. II):\\
A divisor $D$ on $X$ determines a line bundle on $X$, denoted by
$[D]$. The sheaf of sections of $[D]$ may be identified with the
$\mathcal{O}_X$-subsheaf of $\mathcal{M}$ (i.e., a
$\mathcal{O}_X$-submodule of $\mathcal{M}$) generated on the open
cover $U_\iota$ by $1/f_\iota$. Equivalently, transition functions
for $[D]$ with respect to the covering $U_\iota$ are
$g_{\iota\kappa}=f_\iota/f_\kappa$. And a section $\sigma =
\{s_\iota\}, s_\iota = g_{\iota \kappa}s_\kappa$ of a line bundle
$L$ on $X$
determines a divisor $D = \{(U_\iota, s_\iota)\}$.\\

Let $X,Y$ be non-singular compact complex varieties, $f: X \mapsto
Y$ be a proper holomorphic map. Recall that the pushforward
operator on $(p,q)$-currents on $X$, $f_*:\mathcal{D}'^{(p,q)}(X)
\mapsto \mathcal{D}'^{(p,q)}(Y)$ is defined by $<f_*T,\phi> =
<T,f^*\phi>$ for all $T \in \mathcal{D}'^{(p,q)}(X)$ and all $\phi
\in \mathcal{D}^{(p,q)}(Y)$. This operator commutes with the
differential operators $d$ and $d^c$, hence $dd^cf_* = f_*dd^c$.
If $T = dd^cu$, where $u$ is a plurisubharmonic function on $X$
and $f$ has finite fibers, then $f_*u(y) = \sum_{x \in
f^{-1}(y)}u(x), \quad y \in Y$.\\

 For $X,Y$ and $f$ as above it is possible to define (up
to an isomorphism) the pushforward $f_*L$ of a line bundle $L$
over $X$ corresponding to a sheaf $\mathcal{F}$ as the line bundle
corresponding to the sheaf $f_*\mathcal{F}$. The pushforward of an
invertible sheaf is constructed as follows (\cite{GR}, theorem
2.3.4 + Appendix A): Let $y \in Y$ and let $x_1,\ldots x_t$ be
different points in the fiber $f^{-1}(y)$, $U'_1,\ldots,U'_t$-
pairwise disjoint open neighborhoods of $x_1,\ldots,x_t$ and $V'$-
an open neighborhood of $y$. Then there exists $V \subset V'$- an
open neighborhood of $y$ such that $f^{-1}(V) = \bigcup_{j=1}^t
U_j$, where $U_j = f^{-1}(U_j) \cap U'_j$ are pairwise disjoint
neighborhoods of the points $x_j$ , $j=1,\ldots,t$ (in particular,
$f(U_j) \subset V, \quad j=1,\ldots,t$). For a sheaf $\mathcal{F}$
on $X$ these give canonical bijections
$(f_{U_j,V})_*(\mathcal{F}(U_j) \rightleftarrows
\mathcal{F}({f_{U_j,V}}^{-1}(V))=\mathcal{F}(U_j), \quad
j=1,\ldots,t$ (where $f_{U_j,V}= f\mid_{U_j}:U_j \mapsto V$ and
$\mathcal{F}(f^{-1}(V))\rightleftarrows
\prod_{j=1}^t\mathcal{F}(U_j)$ ($\prod$ denotes here the Cartesian
product). These bijections allow us to define $f_*(\mathcal{F})(V)
= \mathcal{F}(f^{-1}(V))$.\\

For pushforward of divisors and associated line bundles, the
following proposition holds:\\
\begin{prop}\label{prop:pushdiv} Let $X,Y$ be non-singular compact complex varieties, $f: X \mapsto
Y$ be a proper holomorphic surjection and let $D$ be a divisor in
$X$. Then:\\
(a)The current $f_*D$ is a divisor in
$Y$;\\
(b) $f_*[D] = [f_*D]$, where $[D]$ is the line bundle associated
with the divisor $D$.
\end{prop}
\begin{proof} Part (a) is Lemma 3.1 from \cite{Ji}. From its proof it can be deduced (cf. also \cite{Fu}, section 1.4) that
if $g_{\iota\kappa}$ are transition functions in $[D]$, then
$f_*(g_{\iota\kappa})$ are transition functions in $[f_*D]$, which
implies part (b).
\end{proof}

\section{Singular hermitian metrics on line bundles and $\omega$- plurisubharmonic functions} \label{s:metric}

In the treatment of singular hermitian metrics on line bundles we
will follow the approach according to Demailly (cf. \cite{De2},
\cite{De3}).\\
\begin{defin}(\cite{De3}, Definition 3.12) Let $L$ be a complex line bundle over a complex
manifold $X$. A singular (hermitian) metric on $L$ is a metric
which is given in any trivialization $\theta: L\mid_{U} \mapsto U
\times \mathbb{C}$ by
\[
\|\xi\| = |\theta(\xi)|\exp(-\psi(x)), \quad x \in U, \xi \in L_x,
\]
where $\psi \in L^1_{loc}(U)$ is a function called the weight of
the metric with respect to the trivialization $\theta$.
\end{defin}
A singular metric can be thus given by a collection of functions
$\psi = \{\psi_{\iota}\}, \quad \psi_{\iota} \in
L^1_{loc}(U_{\iota})$ satisfying $\psi_{\iota} = \psi_{\kappa} +
\log|g_{\iota \kappa}|$ in $U_{\iota \kappa}$, where $U_{\iota}$
is a trivializing cover for $L$ and $g_{\iota \kappa}$ are the
transition functions. The metric is called positive if the
functions $\psi_{\iota}$ are plurisubharmonic.\\
The following proposition is simple but useful:\\
\begin{prop}\label{prop:pushmetric} If the collection
$\{\psi_{\iota}\}$ defines a positive singular metric in a line
bundle $L$ on $X$ and $f:X \mapsto Y$ is a proper holomorphic
surjective map, then $\{f_*(\psi_\iota)\}$ defines a positive
singular metric in the line bundle $f_*L$ on $Y$.
\end{prop}
\begin{proof} If $f$ is holomorphic and proper, the functions $f_*(\psi_\iota)$
are plurisubharmonic if $\psi_\iota$ are (see e.g. \cite{De1},
proposition 1.13). Further, $\psi_\kappa = \psi_\iota + \log
|g_{\iota \kappa}|$ gives $f_*\psi_\kappa(z)= f_*\psi_\iota(z) +
\log \prod_{x \in f^{-1}(z)}|g_{\iota \kappa}(x)|$ in $U_{\iota
\kappa}$. The latter term is the logarithm of the modulus of a
transition function in the bundle $f_*L$ (cf. Proposition
\ref{prop:pushdiv}), so the compatibility conditions also hold.
\end{proof}

Let now $D$ be an effective divisor on $X$ and $[D]$ its
associated line bundle. One can give a positive singular metric on
$[D]$ as in Example (3.14) in \cite{De3}. Assume that
$\sigma_0,\sigma_1,\ldots, \sigma_N$ are nonzero holomorphic
sections of $[D]$. They define a singular hermitian metric on
$[D]$ by
\[
\|\xi\|^2_{\sim} =
\frac{|\theta(\xi)|^2}{|\theta(\sigma_1(x))|^2+\ldots+\theta(\sigma_N(x))|^2}
\] with respect to  a trivialization $\theta$. The weight function
for this metric is
$\psi(x)=\log(\sum_{j=1}^N)|\theta(\sigma_j(x))|^2)^{1/2}$, which
is a plurisubharmonic
function.\\
This metric can be viewed as introduced by a metric on
$H^0(X,[D])$ as follows: Let $\sigma_0,\sigma_1,\ldots,\sigma_N$
be a basis for the linear system $|D|$ of all divisors linearly
equivalent to $D$ and $B_{|D|} = \bigcap {\sigma_j}^{-1}(0)$ its
base locus. There is a meromorphic map $\Phi_{|D|}: X \setminus
B_{|D|} \mapsto \mathbb{P}^{N}, \quad \Phi_{|D|}(x) =
(\sigma_0(x):\ldots:\sigma_N(x))$. Then the positive closed
$(1,1)$-current $\frac{i}{2\pi}\Theta(F) = dd^c \psi$ is equal to
the pullback over $X\setminus B_{|D|}$ of the Fubini-Study metric
$\omega_{FS} = \frac{1}{2\pi}dd^c\log(|z_0|^2+\ldots+|z_N|^2)$ of
$\mathbb{P}^{N}$ by $\Phi_{|D|}$.

In what follows we will assume that $D$ is very ample, i.e., that
the map $\Phi_{|D|}:X \mapsto \mathbb{P}^{N}$ associated with the
linear system $|D| = \mathbb{P}(H^0(X,[D]))$ is a regular
embedding. Then in particular $B_{|D|} = \emptyset$. According to
Theorem II.7.1 in \cite{Ha}, $\Phi^*(\mathcal{O}(1))$ is an
invertible sheaf on $X$, which is generated by the global sections
$\sigma_0,\ldots,\sigma_N$, $\sigma_j=\Phi^*(z_j), \quad
j=0,\ldots,N$. Moreover, as $\Phi$ is an embedding, each open set
$X_j := X \setminus {\sigma_j}^{-1}(0), \quad j=0,\ldots,N$ is
affine (\cite{Ha}, Proposition II.7.2). \\
We will also assume that $X$ is normal, i.e., the linear system
$|D|$ on $X$ giving the embedding $\Phi:X \mapsto \mathbb{P}^N$ is
complete. This means (cf. \cite{GH}, p. 177) the restriction map
$H^0(\mathbb{P}^N,\mathcal{O}(H)) \mapsto H^0(X,\mathcal{O}(H)$ is
surjective, where $H$ denotes the hyperplane bundle.\\

It will be convenient to consider another positive singular metric
on the line bundle $[D]$ associated with a divisor $D$, which can
be introduced using an isomorphism between $H^0(X,[D])$ and $\{u
\in \mathcal{M}(X)^*: \mbox{ div }u \geq -D\}$. Let
$\sigma_0,\sigma_1,...,\sigma_N$ be a basis for $H^0(X,[D])$ We
can assume that $D = \mbox{div }\sigma_0$, so $\mbox{supp }D = \{x
\in X: \sigma_0(x) = 0\}$. Then, for $u = \tau/\sigma_0$, one just
takes $\|u\| = |u| = |\tau|\exp(-\log|\sigma_0|)$ (\cite{De3},
Example 3.13). Hence the weight for this metric coincides with the
logarithm of the local equation of $D$ in the trivialization
$\theta_\iota$ for $[D]$. The following relation holds:\\
\begin{prop}\label{prop:equiv} The metrics $\|\cdot\|_{\sim}$
and $\|\cdot\|$ are equivalent, i.e., $\exists c > 0: \forall \xi:
c^{-1}\|\xi\| \leq \|\xi\|_{\sim} \leq c\|\xi\|$.
\end{prop}
\begin{proof} One uses Rudin- Sadullaev estimates on $X$ in
exactly the same way as in Lemma 2 in \cite{St2}.
\end{proof}

The metric $\|\cdot\|$ is particularly useful in measuring the
distance from a point in $X$ to the support of $D$. Consider the
subspace $Y = \mbox{span}(\sigma_1,\ldots,\sigma_N)$ of
$H^0(X,[D])$. Then $\mathbb{C}\sigma_0 \times \{(0,\ldots,0)\}$ is
the linear complement of $\{0\}\times Y$ in $H^0(X,[D])$. Taking
the norm $|||(\tau\sigma_0,z)||| = |\tau| + \|(0,z)\|$ in
$\mathbb{C}\sigma_0 \times Y$, where $\|\cdot\|$ is the metric in
$H^0(X,[D])$ introduced above, one can construct the following
neighborhood base for the set $Y_{\infty} = \mathbb{P}(0 \times
Y)$, (which equals the support of $D$):
\[
\Omega_K = \{\lambda \in \mathbb{P}(\mathbb{C}\sigma_0\times Y):
\lambda \subset \{(\tau,y): |\tau| < (1/K)\|y\|\}, K>0\},
\]
which is the same as $\{z \in \mathbb{C}^N: |z| \geq K\} \cup
Y_{\infty}$ (\cite{L}, VII.3.4).\\
Let $H$ be a hyperplane in the linear space $Y$ such that $0
\notin H$. There is a unique linear form $\lambda_H \in Y^*$ such
that $H = \{z: \lambda_H = 1\}$. Let$H_* = \ker \lambda_H$. On the
space $\mathbb{P}(\mathbb{C}\times Y$ we have a chart $\beta_H:
\mathbb{P}(\mathbb{C} \times Y)\setminus
\mathbb{P}(\mathbb{C}\times H_*) \ni \mathbb{C}z \mapsto
z/\lambda_H(z) \in \mathbb{C} \times H$. The domains of the charts
$\beta_H$ cover the set $Y_{\infty}=\mathbb{P}(\{0\}\times Y)$.
Having chosen a norm in $Y$ and the norm $|||(t,z)|||=|t|+\|z\|$
in $\mathbb{C}\times H_*$, we can take
$\rho_{\beta_H}(z,Y_\infty\setminus \mathbb{P}(\mathbb{C}\times
H_*) = |\lambda_H(z)|^{-1}, \quad z \in Y \setminus H_*$, as the
distance between $Y_\infty$ and $z \in Y$. In the set $H_i = X_i
\setminus (X_i \cap \{\sigma_o=0\})$,
$\lambda_{H_i}=(\sigma_i/\sigma_0)$ is the reciprocal of
the local equation of $D$, $i=1,\ldots,N$\\

The neighborhoods $\Omega_K$ have an alternative interpretation:\\
\begin{prop}\label{prop:tubes} $ \Omega_K = \{(x,\xi):
|\beta_{H_i}(\xi)|\exp(-\psi(x)) < 1/K\}$ is a strongly
pseudoconvex $1/K$-tube around $A$.
\end{prop}
\begin{proof}
Note that the existence of a strongly pseudoconvex neighborhood of
the zero section is  equivalent to the negativity of a holomorphic
line bundle, see \cite{FG}, Proposition VI.6.1 and VI.6.2. We will
use $\mathcal{O}(-1)$, the tautological line bundle  over the
projective space $\mathbb{P}(\mathbb{C}\times \mathbb{C}^N) =
\mathbb{P}(H^0(X,[D]))$. A point $\lambda \in
\mathbb{P}(\mathbb{C}\times \mathbb{C}^N)$  can be identified with
the fiber of $\mathcal{O}(-1)$.  Under this identification, the
set $Y_{\infty} = A$ corresponds to the zero section of
$\mathcal{O}(-1)$ and $\Omega_K$  is a strongly pseudoconvex
$1/K$-tube around $A$, for all $K > 0$.
\end{proof}

We will also need a notion of $\omega$-plurisubharmonic functions
on a compact connected K\"{a}hler manifold $X$ with respect to a
closed real current $\omega$ of bidegree $(1,1)$ on $X$. ( On
complex projective spaces such functions were first considered in
\cite{BT}).  Under the assumptions we made, $X$ is a projective
algebraic manifold. The class of $\omega$-plurisubharmonic
functions is defined as follows
(see \cite{GZ}):\\
\[
PSH(\omega,X)= \{\phi \in
L^1(X,\mathbb{R}\cup\{-\infty\}):dd^c\phi \geq -\omega, \phi
\mbox{ is upper semicontinuous}\}
\]
Let $L$ be a holomorphic line bundle on $X$ and $h = \{h_\iota\}$
a smooth metric $h = \{h_\iota\}$ in $L$. The curvature $\omega =
dd^ch_\iota$ in $U_\iota$ defines a real closed current globally
on $X$, because of compatibility conditions for $h_\iota$. It was
observed in \cite{GZ} (at the beginning of \S 5) that the class
$PSH(X,\omega)$ is in 1-to-1 correspondence with the set of
positive singular metric in $L$ on $X$ (that is, with the set of
their weights). Namely, for a metric $\psi$ the function $\phi =
\psi - h$ satisfies $dd^c\phi \geq -\omega$, and conversely, for a
$\phi \in PSH(X,\omega)$, the collection $\psi =\{\psi_\iota =
\phi + h_\iota\}$ defines a positive singular metric in $L$ on
$X$.

\section{Attracting divisors}\label{s:attract}
We will use the following definition of an attracting set, which is of topological character:\\
\begin{defin}(see \cite{St1} or \cite{St2} and the references there) Let $(X,d)$
be a metric space and let $f: X \mapsto X$ be a continuous map. A
closed set $A \subset X$ is attracting for $f$ if and only if
$f(A)=A$ and there exists a neighborhood $U$ of $A$, of the form
$U= \{x\in X: d(x,A) \leq \varepsilon\}$, such that $f(U) \subset
\mbox{int }U$ and $\bigcap_{n \geq 0}f^n(U) = A$.
\end{defin}

Recall our working setting: $X$ is a normal projective algebraic
manifold, $D$ is a very ample divisor on $X$ given by a
holomorphic section $\{\sigma_0 =0\}$ (by Kodaira Embedding
Theorem, $D$ is positive), $A$ is the support of $D$, $f:X \mapsto
X$ is a proper holomorphic map with finite fibers, satisfying
$f(A) = A$. Then a neighborhood base for $A$ is given by $\{x:
\psi(x) \leq \varepsilon\}$, where $\psi$ is the weight of a
positive singular metric on the line bundle $[D]$. We will say
that $D$ is attracting for $f$ if its support $A$ is an attracting
set
for $f$ in the sense of the above definition.\\

\begin{thm} The following are equivalent:\\
(i) $D$ is attracting for $f$;\\
(ii) there is a positive singular metric with weight $\psi$ on the
line bundle $[D]$ and a neighborhood $V$ of $A$ such that $\exists
0 < \beta < 1: \quad
f_*\psi \geq \beta \psi$ in $V$;\\
(iii)there is a positive singular metric with weight $\phi$ on the
line bundle $[D]$ such that $\exists 0 < \beta < 1: \quad
f_*(PSH(X,\omega)) \subset \beta PSH(X,\omega)$, where $\omega =
dd^c\phi$.
\end{thm}
\begin{proof}
$(i)\Leftrightarrow (ii)$: Observe that $-\psi(x) =
-\log|\sigma_0(x)| = -\log \mbox{ dist}(x,A)$ is a
plurisubharmonic exhaustion on
$X \setminus A$. The equivalence is the same as Lemma 1 in \cite{St2}.\\
$(ii) \Rightarrow (iii)$: Let $\psi$ be a metric satisfying the
condition (ii) and let $\omega = dd^c\psi$. As in Example 1.2 in
\cite{GZ}, one can establish the following one-to-one
correspondence between $PSH(X,\omega)$ and the class $\{u \in
PSH(X\setminus A): u - \psi \leq c_u\}$: if $u \in PSH(X\setminus
A)$ satisfies $u \leq \psi +c_u$ with some constant $c_u$
depending only on $u$, then $v(x)
=\begin{cases}u(x)-\psi(x),x\notin A\\ \limsup_{(X\setminus A)\ni
y\to x}(u(y)-\psi(y))\end{cases}$ is in $PSH(X,\omega)$. A similar
correspondence can be established for $PSH(X,f_*\omega)$,
involving the weight $f_*\psi$. Let $v \in PSH(X,\omega)$. By
(ii), $f_*v - \psi$ has a well-defined $\limsup_{(X\setminus A)
\ni y \to x}(f_*v(y)-\psi(y)) \leq \limsup_{(X\setminus A)\ni y
\to x}(f_*v(y) -f_*\psi(y))$ for $x \in A$, hence $f_*v \in
PSH(X,\omega)$.\\
$(iii) \Rightarrow (ii)$: By (iii), $f_*(-\phi) \in PSH(X,\beta
\omega)$, so $f_*(-\phi) + \beta\phi$ defines a metric in $[D]$.
On the set supp $D$ this metric is equal to $0$, hence $f_*\phi -
\beta\phi \to \infty$ as $x \to A$, which implies that close
enough to $A = \mbox{supp }D$, $f_*\phi - \beta\phi > 0$. Thus
$\psi
=\phi$ is a metric as in condition (ii).\\
\end{proof}

\begin{rk}(1) Even though trivializations were used to formulate the definition of a singular
metrics, the Theorem is independent of trivializations defining
the line bundle $L = [D]$. Indeed, two sets of trivializations
$\theta_\iota, \theta'_\iota$ defining $L$ in the open cover
$U_\iota$ are related by $\theta'_\iota = \eta_\iota \cdot
\theta_\iota$, with $\eta_\iota \in \mathcal{O}^*(U_\iota)$, so
the form $\omega$ in the Theorem will be replaced by a
cohomologous form $\omega' = \omega + dd^c\chi$, where $\chi$ is a
function on $X$ which is integrable with respect to any smooth
volume form. Accordingly, the metric $\psi$ will be replaced by
$\psi - \chi$ and the class $PSH(X,\omega)$ by $PSH(X,\omega')=
PSH(X,\omega) + \chi$ (the last equality is Proposition 1.3.3 in
\cite{GZ}).\\
(2) It was observed in \cite{DS}, section 6, that the classes of
$\omega$-plurisubharmonic functions are in general not stable
under the operator $f_*$. Condition (iii) of our Theorem gives an
example of situation in which stability occurs, since $\beta PSH(X,\omega)\subset PSH(X,\omega)$ for $0<\beta<1$.\\
(3) If $A$ is a regular hypersurface in $X$, then the first
adjunction formula (see e.g. \cite{FG}, IV.5.10, \cite{GH}, p.
145) says that on $A$ the line bundle $[A]$ coincides with the
normal bundle $N_X(A)$ of $A$ in $X$. The neighborhoods in the
definition of the attracting set and in the Theorem can then be
chosen as tubular neighborhoods of $A$ (cf. \ref{prop:tubes}.\\
\end{rk}

\end{document}